\documentclass[12pt,a4paper,fullpage]{article}
\usepackage{amsmath,amsthm}
\usepackage{amssymb}

\newtheorem{theorem}{Theorem}
\newtheorem{corollary}{Corollary}

\newtheorem{remark}{Remark}

\title{{\bf A class number problem for imaginary cyclic number fields of 2-power degrees}}
\author{
St\'ephane R. LOUBOUTIN\\
Aix Marseille Universit\'e, CNRS, Centrale Marseille, I2M,\\ 
Marseille, France\\
stephane.louboutin@univ-amu.fr}
\date{\today}

\begin{document}
\bibliographystyle{alpha}
\maketitle

\begin{abstract}
In 2024, M. K. Ram proved that the class number of an imaginary cyclic quartic number field 
is never equal to a prime $p\equiv 3\pmod 4$. 
Here we greatly generalize this result to the case of the non-quadratic imaginary cyclic number fields of $2$-power degrees 
and not necessarily prime class numbers.
\end{abstract}

\section{Introduction}
In 2024, M. K. Ram proved that the class number of an imaginary cyclic quartic number field 
is never equal to a prime $p\equiv 3\pmod 4$, see \cite{Ram}. 
Here in Theorem \ref{mainresult} and Corollary \ref{Lou}, 
we greatly generalize this result to the case of the non-quadratic imaginary cyclic number fields of $2$-power degrees 
and not necessarily prime class numbers. 
For the reminders that follow, we mainly refer to L. C. Washington's book \cite{Was}.

Let $N$ be an imaginary abelian normal field of conductor $f_N>2$. 
Its degree $2n$ is even 
and its subfield $N^+$ fixed by the complex conjugation is totally real. 
The class number $h_{N^+}$ of $N^+$ divides the class number $h_N$ of $N$ 
and the positive integer $h_N^- =h_N/h_{N^+}$ is the relative class number of $N$. 
It divides $h_N$. 
Let $Q_N\in\{1,2\}$ be the Hasse unit index of $N$ 
and $w_N$ be the number of complex roots of unity in $N$.
Let $X_N$ be the group of order $2n$ of primitive Dirichlet characters associated with $N$ 
and let $X_N^-$ denote the subset of cardinality $n$ of the odd characters in $X_N$, 
i.e. $X_N^-=\{\chi\in X_N;\ \chi(-1)=-1\}$. 
Let $L(s,\chi) =\sum_{k\geq 1}\chi(k)k^{-s}$ be the $L$-series associated with $\chi\in X_N^-$. 
It convergences absolutely in the half-open plane $\Re (s)>1$ 
and admits a holomorphic continuation to the whole complex plane.
Then 
\begin{equation}\label{formulahNminus}
h_N^-
=Q_Nw_N2^{-n}\prod_{\chi\in X_N^-}L(0,\chi)
\end{equation} 
(see \cite[(17]{LouNagoya161}).
Now, let $\chi$ be a Dirichlet character of order $n_\chi>1$ and conductor $f_\chi>1$. 
Then $L(0,\chi)\in {\mathbb Q}(\zeta_{n_\chi}^{})$. 
The Galois group ${\rm Gal}({\mathbb Q}(\zeta_{f_\chi}^{})/{\mathbb Q})$ 
is canonically isomorphic to the multiplicative group $({\mathbb Z}/f_\chi {\mathbb Z})^*$ 
by 
$c\in ({\mathbb Z}/f_\chi {\mathbb Z})^*$ 
maps to $\sigma_c\in {\rm Gal}({\mathbb Q}(\zeta_{f_\chi}^{})/{\mathbb Q})$ 
defined by $\sigma_c(\zeta_{f_\chi}^{}) =\zeta_{f_\chi}^{c}$. 
Let $M_\chi$ be the subfield of ${\mathbb Q}(\zeta_{f_\chi}^{})$ fixed by $\ker\chi$, 
i.e. $\alpha \in M_\chi$ if and only if $\alpha\in {\mathbb Q}(\zeta_{f_\chi}^{})$ 
and $\sigma_c(\alpha)=\alpha$ for $c\in ({\mathbb Z}/f_\chi {\mathbb Z})^*$. 
Then $M_\chi$ is of degree $n_\chi$.
Finally let $w_\chi$ denote the number of complex roots of unity in $M_\chi$. 
We have $w_\chi L(0,\chi)\in {\mathbb Z}[\zeta_{n_\chi}^{}]$ 
(see \cite[page 170]{LouNagoya161}).

From now on we assume that $N$ is moreover cyclic and of degree $$2n=2^m\geq 4,$$ 
a perfect power of $2$. 
Then $X_N$ is cyclic of order $2n$ . 
Let $\chi_N$ be any generator of $X_N$ then 
$X_N^- =\{\chi_N^k;\ 1\leq k\leq 2n\text{ and $k$ odd}\}.$
Then as explained above and by \cite[(18)]{LouNagoya161} we have
\begin{equation}\label{formulahNminuscyclic2powerdegree}
w_NL(0,\chi_N)\in {\mathbb Z}[\zeta_{2n}]
\text{ and }
h_N^-
=\frac{Q_Nw_N}{2^n}N_{{\mathbb Q}(\zeta_{2n})/{\mathbb Q}}(L(0,\chi_N)).
\end{equation} 
Finally, by \cite[Lemmas (b) and (c)]{LouMathComp64}, we have
\begin{equation}\label{Qw}
Q_N=1
\text{ and } w_N
=\begin{cases}
2\ell&\hbox{if $2n+1=2^m+1=\ell$ is prime and $N ={\mathbb Q}(\zeta_\ell)$,}\\
2&\hbox{otherwise.}
\end{cases}
\end{equation}
Notice that if $l=2^m+1$ is prime then $m=2^{m'}$ is a perfect power of $2$ 
and $(m,l)\in\{(1,3),(2,5),(4,17),(8,257),(16,65537),\cdots\}$.

\section{There are only finitely many imaginary cyclic numbers fields of 2-power degrees 
with some restrictions on their class numbers}
Here is our main result:

\noindent\frame{\vbox{
\begin{theorem}\label{mainresult}
Let $N$ be an imaginary cyclic number field of $2$-power degree $2n=2^m\geq 4$ 
whose class number $h_N$ is of the form $h_N=2^ad$, 
where $d$ is odd, square-free and every prime divisor $p$ of $d$ satisfies $p\not\equiv 1\pmod {2n}$.
Then $h_N^-$ divides $2^a$.
Consequently, for a given $a\geq 0$ there are only finitely many such imaginary cyclic number fields 
and we can compute an explicit upper bound $B_a$ on their conductors. 
For example, 
$B_0=2 500$, 
$B_1 =6 300$, 
$B_2 =16 000$, 
$B_3 =36 000$
and $B_4 =84 000$.
\end{theorem}
}}

\begin{proof}
First, let $p\geq 3$ be an odd prime. 
In the cyclotomic number field ${\mathbb Q}(\zeta_{2n})$ 
it splits into $n/f$ distinct prime ideals 
each of residue degree $f$ and norm $p^f$, 
where $f$ is the order of $p$ in the multiplicative group $({\mathbb Z}/2n{\mathbb Z})^*$ 
(see \cite[Theorem 2.3]{Was}). Hence, $f=1$ if and only if $p\equiv 1\pmod {2n}$, 
i.e. the norm of a prime ideal ${\mathcal P}$ of ${\mathbb Q}(\zeta_{2n})$ above an odd prime $p$ 
is square-free if and only if $p\equiv 1\pmod {2n}$. 

Second, take $\alpha \in {\mathbb Z}[\zeta_{2n}]$.
By looking at the prime ideal factorization $(\alpha) ={\mathcal P}_1^{e_1}\cdots {\mathcal P}_r^{e_r}$ 
of the principal ideal $(\alpha) =\alpha {\mathbb Z}[\zeta_{2n}]$,
we deduce that if an odd prime $p\not\equiv 1\pmod{2n}$ divides the norm 
$N_{{\mathbb Q}(\zeta_{2n})/{\mathbb Q}}(\alpha)$ of $\alpha \in {\mathbb Z}[\zeta_{2n}]$
then $p^2$ divides it. 

Third, let $N$ be an imaginary cyclic number field of $2$-power degree $2n=2^m\geq 4$. 
By \eqref{formulahNminuscyclic2powerdegree} and \eqref{Qw}, 
if $2n+1=2^m+1=\ell$ is prime and $N ={\mathbb Q}(\zeta_\ell)$ we have
\begin{equation}\label{2n+1=l}
2^{2n-1}\ell^{n-1}h_N^-
=N_{{\mathbb Q}(\zeta_{2n})/{\mathbb Q}}(2\ell L(0,\chi_N))
\text{ and }
2\ell L(0,\chi_N)\in {\mathbb Z}[\zeta_{2n}]
\end{equation} 
and otherwise we have
$$2^{2n-1}h_N^-
=N_{{\mathbb Q}(\zeta_{2n})/{\mathbb Q}}(2L(0,\chi_N))
\text{ and }
2L(0,\chi_N)\in {\mathbb Z}[\zeta_{2n}].$$
Noticing that in case \eqref{2n+1=l} we have $\ell\equiv 1\pmod {2n}$, 
it follows that if an odd prime $p\not\equiv 1\pmod{2n}$ divides $h_N^-$, 
then its square $p^2$ divides $h_N^-$. 

Now, the first assertion follows. 
Indeed if $h_N^-$ that divides $h_N$ were not a power of $2$, 
there would exist a prime $p\geq 3$ dividing $h_N^-$. 
Then $p$ would divide $h_N$, 
hence would divide $d$. 
We would have $p\not\equiv 1\pmod {2n}$ 
and consequently $p^2$ would divide $h_N^-$. 
Therefore, $p^2$ would divide $h_N$, 
hence would divide $d$, 
contradicting that $d$ is square-free.

For the second assertion, we use \cite[Corollary 20]{LouActa121} 
and notice that $f_N\geq 2500$ implies $\rho_N\geq\sqrt{f_N}\geq 50$, 
by \cite[Corollary 1]{Mur}.
\end{proof}

\begin{remark}
The referee came up with the following more algebraic proof of the main step of Theorem \ref{mainresult}. 
Suppose $p\not\equiv 1\pmod{2^m}$ divides $h_N$ exactly to the first power. 
Let $A_p \backsimeq {\mathbb Z}/p{\mathbb Z}$ 
 be the Sylow $p$-subgroup of the class group. 
 By assumption, $G = {\rm Gal}(N/{\mathbb Q})$
is cyclic of order $2^m$ and has $J =$ complex conjugation as its unique element of order $2$. 
The group $G$ acts on $A_p$ in the usual way. 
Since the automorphism group of $A_p$ has no elements of order $2^m$, 
we must have some non-identity elements of $G$ act trivially. 
In particular, $J$ acts trivially. 
Therefore, $1+J$ acts by squaring. 
But $1 + J$ annihilates the minus part of the class group, 
so squaring annihilates $A_p^-$. 
Since $p$ is odd, this implies that $A_p^- = 1$.
\end{remark}

As a Corollary we obtain a great generalization of Ram's result in \cite{Ram} 
which asserts that the class number $h_N$ of an imaginary cyclic quartic number field 
is never a prime numbers $p\equiv 3\pmod 4$:

\noindent\frame{\vbox{
\begin{corollary}\label{Lou}
The class number $h_N$ of an imaginary cyclic number field of $2$-power degree $2n=2^m\geq 4$ 
is never a square-free integer product of prime numbers $p\not\equiv 1\pmod {2n}$, 
and never twice, four times, eight times or sixteen times such a product.
\end{corollary}
}}

\begin{proof}
The results follows from Theorem \ref{mainresult}, 
according to which $h_N^-=1$, $2$, $4$, $8$ or $16$ 
if the class number of $N$ is as required in the statement of this Corollary, 
and to the Tables given in \cite{PK}, 
according to which $h_N^+$ is equal to $1$, $2$ or $4$ whenever $h_N^-=1$, $2$, $4$, $8$ or $16$.
\end{proof}

\bibliography{central}

\end{document}